\documentclass{amsart} % Nicer than default article style: less
\usepackage{amsmath,amsthm}     % Handy math stuff, theorem environments.
\usepackage{amssymb}            % Fancy math symbols.
\usepackage{euscript}           % Nice script font.
\usepackage{enumerate,calc}
\usepackage[matrix,arrow,curve,frame]{xy}    % XY-pic diagram pac

\xymatrixcolsep{1.9pc}                          % Adjust size of diagrams.
\xymatrixrowsep{1.9pc}
\newdir{ >}{{}*!/-5pt/\dir{>}}                  % Make better tailed arrows

\newtheorem{thm}[subsection]{Theorem}
\newtheorem{defn}[subsection]{Definition}
\newtheorem{prop}[subsection]{Proposition}

\newtheorem{cor}[subsection]{Corollary}
\newtheorem{lemma}[subsection]{Lemma}

\newtheorem{ques}[subsection]{Question}

\theoremstyle{definition}  % Bold headings and Roman body text.

\newtheorem{remark}[subsection]{Remark}

  {\end{list}}
  
\newcommand{\dfn}{\textbf} % Make defined words bold.

\newcommand{\mdfn}[1]{\dfn{\mathversion{bold}#1}} % Even make math bold

\newcommand{\field}[1]  {\mathbb #1} % Use blackboard bold for these sets

\newcommand{\R}         {\field R}

\newcommand{\HH}        {\field H}
\newcommand{\OO}        {\field O}

\newcommand{\C}         {\field C}

\DeclareMathOperator{\Ann}{Ann}

\DeclareMathOperator{\Eig}{Eig}

\newcommand{\norm}[1]{\mid \! #1 \! \mid}       %\norm{x} gives |x|

\newcommand{\rea}[1]{|{#1}|}             %geometric realization of #1
\newcommand{\map}{\rightarrow}

\newcommand{\ceck}[1]{\Cech(#1)}         %Cech complex for #1
\newcommand{\oceck}[1]{\Cech^{o}(#1)}    %Ordered Cech complex for #1
\newcommand{\oreal}[1]{\rea{\oceck{U}}}  %Realization of ordered Cech cplex
\newcommand{\creal}[1]{\rea{\ceck{U}}}   %Realization of the Cech complex
\newcommand{\Cech}{\check{C}}

\renewcommand{\Re}{\text{Re}}
\renewcommand{\Im}{\text{Im}}

\newcommand{\llangle}{\langle\langle}
\newcommand{\rrangle}{\rangle\rangle}

\numberwithin{equation}{subsection}

\newenvironment{myequation}
  {\addtocounter{subsection}{1}\begin{eqnarray}}
  {\end{eqnarray}$\!\!$}

\begin{document}

\title{Large annihilators in Cayley-Dickson algebras II}

\author[D.~K.~Biss]
{Daniel K. Biss}\thanks{This research was conducted during the
period the first author served as a Clay Mathematics Institute
Research Fellow}
\author[J.~D.~Christensen]{J.~Daniel Christensen}
\thanks{The second author was partially supported by NSERC}
\author[D.~Dugger]{Daniel Dugger}
\thanks{The third author was partially supported by the National Science Foundation}
\author[D.~C.~Isaksen]{Daniel C. Isaksen}
\thanks{The fourth author was partially supported by the National
Science Foundation.}

\address{Department of Mathematics\\
%5734 S. University Ave.\\
University of Chicago\\Chicago, IL 60637}

\address{Department of Mathematics\\ 
University of Western Ontario\\
London, Ontario N6A 5B7\\
Canada}

\address{Department of Mathematics\\ University of Oregon\\ Eugene, OR
97403}

\address{Department of Mathematics\\ Wayne State University\\
Detroit, MI 48202}

\email{daniel@math.uchicago.edu}

\email{jdc@uwo.ca}

\email{ddugger@math.uoregon.edu}

\email{isaksen@math.wayne.edu}

\begin{abstract}
We establish many previously unknown properties of zero-divisors in
Cayley-Dickson algebras.  The basic approach is to use a certain
splitting that simplifies computations surprisingly.
\end{abstract}

\maketitle

\section{Introduction}

Cayley-Dickson algebras are non-associative finite-dimensional
$\R$-division algebras that generalize the real numbers, the complex
numbers, the quaternions, and the octonions.  This paper is a sequel
to \cite{DDD}, which explores
some detailed algebraic properties of these algebras.

Classically, the first four Cayley-Dickson algebras, i.e.,
$\R$, $\C$, $\HH$, and $\OO$, are viewed as at least somewhat
well-behaved, while
the larger Cayley-Dickson algebras are considered to be pathological.
There are several different ways of making this distinction.  One
difference is that the first four algebras do not possess zero-divisors,
while the higher algebras do have zero-divisors.  Our primary
long-term goal is to understand the zero-divisors in as much detail
as possible.  The specific purpose of this paper is to build 
directly on the ideas of \cite{DDD} about zero-divisors with
large annihilators.
Our motivation for studying zero-divisors
is their potential for useful applications in topology; see
\cite{Co} for more details.

Let $A_n$ be the Cayley-Dickson algebra of dimension $2^n$.
The central idea of the paper is to use a certain additive
splitting of $A_n$ (as expressed indirectly in Definition \ref{defn:bracket})
to simplify multiplication formulas.  
Multiplication does not quite respect the splitting, but it almost does
(see Proposition \ref{prop:bracket-multiply}).
Theorem \ref{thm:bracket-multiply} is the technical heart of the paper; it
supplies expressions for
multiplication of elements of a codimension 4 subspace of $A_n$
that are simpler than one might expect.

These simple multiplication formulas lead to detailed information about
zero-divisors and their annihilators.
Section \ref{sctn:ann-H-perp} takes a straightforward approach:
just write out equations and solve them as explicitly
as possible.  Our simple multiplication formulas make this feasible.
This leads to 
Theorem \ref{thm:ann-bracket-bound}, which almost completely
computes the dimension of the annihilator of any element.
There are two ways in which the theorem
fails to be complete.  
First, it only treats annihilators of 
elements in a codimension 4 subspace of $A_n$.
Second, rather than determining the dimension of an annihilator
precisely, it gives two options, which differ by 4.

We currently have no solution to the first problem.
However, in this regard, it was already known that one codimension 2 slice
is easy to deal with, so the restriction is really only codimension 2.
We intend to address this question in future work.

The second problem has a partial solution in 
Theorems \ref{thm:ann-off-D-locus} and \ref{thm:ann-D-locus},
which distinguish between the two possible cases.
We find that the answer for $A_{n+1}$ depends inductively 
not just on an understanding of zero-divisors in $A_n$ but also on a
detailed understanding of annihilators in $A_n$
(see Definition \ref{defn:D}).
Therefore, the description in these theorems is not as explicit as we
might like.  

Fortunately, we have a complete understanding of zero-divisors and
their annihilators
in $A_4$ \cite[Sections 11 and 12]{DDD}.
This allows us to make calculations about
zero-divisors in $A_5$ that are not yet possible for $A_n$ with $n \geq 6$.
Section \ref{sctn:D5} contains the details of these calculations
in $A_5$.
Consequently, even though we have not made this result explicit in this
article, it is possible to completely understand in geometric terms
the zero-divisors in a codimension 4 subspace of $A_5$.  
This goes
a long way towards completely describing the zero-divisors of $A_5$.

In addition to the concrete results in Section \ref{sctn:D5} about $A_5$,
Section \ref{sctn:stable}
gives a number of results about spaces of zero-divisors in
$A_n$ for arbitrary $n$.  Consider for a moment only 
the zero-divisors whose annihilators have dimension differing 
from the maximum possible dimension by a fixed constant.
We show in Theorem \ref{thm:stable-dim} that, in a certain sense,
the space of such zero-divisors does not depend on $n$.
This is a kind of stability result for zero-divisors with large
annihilators; it was alluded to in~\cite[Remark 15.8]{DDD}.
The basic approach is to use the previous calculations of dimensions
of annihilators, together with bounds on the dimensions of annihilators
from \cite{DDD} (see Theorem \ref{thm:4dim}).

The paper contains a review in Section \ref{sctn:cd} of 
the key properties of Cayley-Dickson algebras that we will use.
Only some of the material is original; it quotes many
results from \cite{DDD} that will be relevant here.

We make one further remark about generalities.  Many of our results
have hypotheses that eliminate consideration of the classical
algebras $\R$, $\C$, $\HH$, and $\OO$, even though sometimes this is not
strictly necessary.  From the perspective of
this paper, these low-dimensional algebras behave 
significantly differently than $A_n$
for $n \geq 4$.  We eliminate them to avoid awkward but easy
special cases.

\subsection{Statement of Results}

We now present a summary of our technical results.

Recall that $A_n$ is additively isomorphic to $A_{n-1} \times A_{n-1}$,
so elements of $A_n$ are expressions $(a,b)$, where $a$ and $b$ belong
to $A_{n-1}$.  See Section \ref{sctn:cd} for a multiplication formula
with respect to this notation.

The element $i_n = (0,1)$ of $A_n$ has many special properties that 
will be described below.  Let $\C_n$ be the linear span of $1 = (1,0)$
and $i_n$; it is a subalgebra of $A_n$ isomorphic to the complex
numbers.  Let $\HH_{n+1}$ be the linear span of 
$(1,0)$, $(0,1)$, $(i_n,0)$, and $(0,i_n)$; it is a subalgebra
of $A_{n+1}$ isomorphic to the quaternions.

It turns out that $A_n$ is naturally a Hermitian inner product space.
The Hermitian inner product of two elements $a$ and $b$ is
the orthogonal projection of $ab^*$ onto $\C_n$.
We say that two elements $a$ and $b$ are $\C$-orthogonal if their
Hermitian inner product is zero.

Results of \cite{DDD} suggest that we should pay particular attention
to elements of $A_{n+1}$ of the form $(a, \pm i_n a)$ with $a$ in the
orthogonal complement $\C_n^\perp$ of $\C_n$.  Every element
of the orthogonal complement $\HH_{n+1}^\perp$ of $\HH_{n+1}$ can be 
written uniquely in the form
\[
\frac{1}{\sqrt{2}} \Big( a, -i_n a \Big) + 
\frac{1}{\sqrt{2}} \Big( b, i_n b \Big),
\]
where $a$ and $b$ belong to $\C_n^\perp$.
We use the notation $\{a,b\}$ for this expression.
We insert the ungainly scalars $\frac{1}{\sqrt{2}}$ in order to
properly normalize some
formulas that appear later.
We would like to consider the product of two elements
$\{a,b\}$ and $\{x,y\}$ of $\HH_{n+1}^\perp$.  

\begin{prop}
\label{prop:intro-mult}
Let $a$, $b$, $x$, and $y$ belong to $\C_n^\perp$, and suppose that
$a$ and $b$ are $\C$-orthogonal to both $x$ and $y$.  Then
\[
\{a, b\} \{x, y \} = \sqrt{2} \{ ax, by \}.
\]
\end{prop}

This result is proved at the beginning of Section \ref{sctn:mult}.
The formula is remarkably simple, but it is not completely general
because of the orthogonality assumptions on $a$, $b$, $x$, and $y$.
Most of Section \ref{sctn:mult} is dedicated to generalizing this 
formula and understanding the resulting error terms.

Recall that the annihilator $\Ann(x)$ of an element $x$ of $A_n$
is the set of all elements $y$ such that $xy = 0$.
Proposition \ref{prop:intro-mult} is the key computational step in
the following theorem about annihilators, which is proved
in Section \ref{sctn:ann-H-perp}.

\begin{thm}
\label{thm:intro-dim-ann}
Let $n \geq 3$, and 
let $a$ and $b$ be non-zero elements of $\C_{n}^\perp$.  Then 
the dimension of the annihilator of $\{a,b\}$ is equal to
$\dim \Ann a + \dim \Ann b$ or 
$\dim \Ann a + \dim \Ann b +4$.
\end{thm}

In order to distinguish between the two cases of 
Theorem \ref{thm:intro-dim-ann}, we need the following definition.

\begin{defn}
\label{defn:intro-D}
The $D$-locus is the space of all elements $\{a,b\}$ of $A_{n+1}$
with $a$ and $b$ in $\C_n^\perp$
such that 
\begin{enumerate}
\item
$a$ and $b$ are $\C$-orthogonal,
\item
$a$ and $\Ann(b)$ are orthogonal, and
\item
$b$ and $\Ann(a)$ are orthogonal.
\end{enumerate}
\end{defn}

The following result is proved in Section \ref{sctn:D-locus}.

\begin{thm}
\label{thm:intro-D-locus}
Let $a$ and $b$ be non-zero elements of $\C_n^\perp$.
If $\{a,b\}$ does not belong to the $D$-locus in $A_{n+1}$, 
then the dimension of the annihilator of $\{a,b\}$ is 
$\dim \Ann a +\dim \Ann b$.
If $\{a,b\}$ belongs to the $D$-locus in $A_{n+1}$, 
then the dimension of the annihilator of $\{a,b\}$ is 
$\dim \Ann a +\dim \Ann b + 4$.
\end{thm}

For example, if neither $a$ nor $b$ are zero-divisors in $A_n$ 
and $a$ and $b$ are not $\C$-orthogonal,
then $\{a,b\}$ is not a zero-divisor.  
If neither $a$ nor $b$ are zero-divisors in $A_n$ but are
$\C$-orthogonal, then $\{a,b\}$ does belong to the $D$-locus in $A_{n+1}$
and is thus a zero-divisor.
On the other hand, the theorem also shows that if $a$ or $b$ is a zero-divisor,
then $\{a,b\}$ is a zero-divisor regardless of whether or not it belongs
to the $D$-locus.  In summary,
if $a$ and $b$ 
are $\C$-orthogonal, then $\{a,b\}$ is always a zero-divisor.

In Section \ref{sctn:D5}, we explicitly work out the meaning of
Definition \ref{defn:intro-D} when $a$ and $b$ belong 
to $A_4$.  The only difficult case occurs when both $a$ and $b$
have non-trivial annihilators, i.e., when both $a$ and $b$ are 
zero-divisors in $A_4$.  This case is explicitly handled in
Theorem \ref{thm:D5}.

Finally, Section \ref{sctn:stable} provides some general results about
zero-divisors with very large annihilators.  Recall from
\cite{DDD} that the largest annihilators in $A_n$ are
$(2^n - 4n + 4)$-dimensional.

\begin{defn}
\label{defn:intro-T}
Let $n \geq 4$, and let $c$ be a multiple of $4$
such that $0 \leq c \leq 2^n - 4n$.
The space \mdfn{$T^c_n$} is the space of elements of length one in $A_n$
whose annihilators have dimension at least
$(2^n - 4n + 4) - c$.  
\end{defn}

In other words, $T^c_n$ consists of the zero-divisors 
with annihilators whose dimensions are within $c$ of the maximum.

\begin{thm}
\label{thm:intro-stable}
Let $n \geq 4$, and let $c$ be a multiple of $4$
such that $0 \leq c \leq 2^n - 4n$.  If $n \geq \frac{c}{4} + 4$,
then $T_{n+1}^c$ is equal to 
the space of elements of the form $\{a,0\}$ or $\{0,a\}$ such that
$a$ belongs to $T_n^c$.
\end{thm}

Theorem \ref{thm:intro-stable}, which is proved in
Section \ref{sctn:stable}, tells us that 
for sufficiently large $n$, the space $T^c_{n+1}$ is diffeomorphic
to the disjoint union of two copies of $T^c_n$.  
The case $c=0$ was proved in \cite[Theorem 15.7]{DDD}.
An interesting
open question is to determine explicitly the geometry of a connected
component of $T^c_n$ for $n$ sufficiently large; this connected
component depends only on $c$.

%%%%%%%%%%%%%%%%%%%%%%%%%%%%%%%%%%%%%%%%%%%%%%%%%%%%%%

\section{Cayley-Dickson algebras}
\label{sctn:cd}

The \mdfn{Cayley-Dickson algebras} are a sequence of non-associative
$\R$-algebras with involution.  
See \cite{DDD} for a full explanation of their basic properties.

These algebras are defined inductively.
We start by defining \mdfn{$A_0$} to be $\R$ with trivial conjugation.  
Given $A_{n-1}$, the algebra
\mdfn{$A_n$} is defined additively to be $A_{n-1} \times A_{n-1}$.
Conjugation in $A_n$ is defined by
\[
(a,b)^* = (a^*,-b),
\]
and multiplication is defined by
\[
(a,b)(c,d) = (ac - d^*b, da + bc^*).
\]

One can verify directly from the definitions that $A_1$ is isomorphic
to the complex numbers $\C$; 
$A_2$ is isomorphic to the quaternions $\HH$; 
and $A_3$ is isomorphic to the octonions $\OO$.  

We implicitly view $A_{n-1}$ as the subalgebra $A_{n-1} \times 0$ of $A_n$.

%%%%%%%%%%%%%%%%%%%%%%%%%%%%%%%%%%%%%%%%%%%%%%%%%%%%%%%%%%%

\subsection{Complex structure}

The element \mdfn{$i_n$} $ =(0,1)$ of $A_n$ 
enjoys many special properties.  One of the primary themes of our long-term
project is to fully exploit these special properties.

Let \mdfn{$\C_n$} be 
the $\R$-linear span of $1 = (1,0)$ and $i_n$.  It is a subalgebra of $A_n$
isomorphic to $\C$.

\begin{lemma}[DDD, Proposition 5.3]
\label{lem:C-vs}
Under left multiplication, $A_n$ is a $\C_n$-vector space.  In particular,
if $\alpha$ and $\beta$ belong to $\C_n$ and $x$ belongs to $A_n$, then
$\alpha(\beta x) = (\alpha \beta) x$.
\end{lemma}

As a consequence, the expression $\alpha \beta x$ is unambiguous;
we will usually simplify notation in this way.

The \mdfn{real part $\Re(x)$} of an element $x$ of $A_n$ is defined to
be $\frac{1}{2}(x + x^*)$, while the 
\mdfn{imaginary part $\Im(x)$} is defined to be $x - \Re(x)$.

The algebra $A_n$ becomes a positive-definite real inner product space
when we define \mdfn{$\langle a, b \rangle_{\R}$} $= \Re(ab^*)$
\cite[Proposition 3.2]{DDD}.
If $a$ and $b$ are imaginary and orthogonal, then $ab$ is imaginary.
Hence, $ba = b^* a^* = (ab)^* = -ab$.  In other words, orthogonal
imaginary elements anti-commute.
A simple calculation shows that $aa^*$ and $a^*a$ are both
equal to $\langle a, a \rangle_{\R}$ for all $a$ in $A_n$
\cite[Lemma~3.6]{DDD}.

We will need the following slightly technical result.

\begin{lemma}
\label{lem:ortho1}
Let $x$ and $y$ be elements of $A_n$ such that $y$ is imaginary.  Then 
$x$ and $xy$ are orthogonal.
\end{lemma}

\begin{proof}
We wish to show that $\Re(x(xy)^*)$ equals zero.
This equals $-\Re( (xx^*) y )$
by \cite[Lemmas~2.6~and~2.8]{DDD}, which
is zero because $y$ is imaginary
and because $xx^*$ is real.
\end{proof}

The real inner product allows us to define a positive-definite 
Hermitian inner product on $A_n$ by setting \mdfn{$\langle a, b \rangle_{\C}$}
to be the orthogonal projection of $ab^*$ onto the subspace $\C_n$ of $A_n$
\cite[Proposition 6.3]{DDD}.
We say that two elements $a$ and $b$ are \mdfn{$\C$-orthogonal}
if $\langle a, b \rangle_{\C} = 0$.

We will frequently consider the subspace
\mdfn{$\C_n^\perp$} of $A_n$; it is the orthogonal complement of $\C_n$
(with respect either to the real or to the Hermitian inner product).
Note that $\C_n^\perp$ is a $\C_n$-vector space; in other words,
if $a$ belongs to $\C_n^\perp$ and $\alpha$ belongs to $\C_n$,
then $\alpha a$ also belongs to $\C_n^\perp$ \cite[Lemma 3.8]{DDD}.

\begin{lemma}[DDD, Lemmas 6.4 and 6.5]
\label{lem:C-conj-linear}
If $a$ belongs to $\C_n^\perp$, then 
left multiplication by $a$ 
is $\C_n$-conjugate-linear
in the sense that $a \cdot \alpha x = \alpha^* \cdot ax$ for any 
$x$ in $A_n$ and any $\alpha$ in $\C_n$.
Moreover, left multiplication is anti-Hermitian in the sense that
$\langle a x, y \rangle_{\C} = -\langle x, a y \rangle_{\C}^*$.
\end{lemma}

Similar results hold for right multiplication by $a$.  
See also \cite[Lemma 2.1]{M2} for
a different version of the claim about conjugate-linearity.

The conjugate-linearity of left and right multiplication
is fundamental to many later calculations.
To emphasize this point, we provide some computational consequences.
The next lemma can be interpreted as a restricted kind of
bi-conjugate-linearity for multiplication.

\begin{lemma}
\label{lem:C-bi-conj}
Let $a$ and $b$ be $\C$-orthogonal elements of $\C_n^\perp$,
and let $\alpha$ and $\beta$ belong to $\C_n$.  Then
$\alpha a \cdot \beta b = \alpha^* \beta^* \cdot ab$.
\end{lemma}

\begin{proof}
By left conjugate-linearity, 
$\alpha a \cdot \beta b = \beta^* (\alpha a \cdot b)$.
Use right conjugate-linearity twice to compute that 
$\beta^* (\alpha a \cdot b) = \beta^* (ab \cdot \alpha)$.
Because $a$ and $b$ are $\C$-orthogonal, $ab$ belongs to $\C_n^\perp$.
Therefore, 
$\beta^* (ab \cdot \alpha) = \beta^* (\alpha^* \cdot ab)$ by
left conjugate-linearity again.
Finally, this equals $\alpha^* \beta^* \cdot ab$ by
Lemma \ref{lem:C-vs}.
\end{proof}

Norms of elements in $A_n$ are defined with respect to either the
real or Hermitian inner
product: $\norm{a} = \sqrt{\langle a, a\rangle_{\R}} =
\sqrt{\langle a, a\rangle_{\C}} = \sqrt{ aa^*}$; this makes sense
because $aa^*$ is always a non-negative real number \cite[Lemma 3.6]{DDD}.
Note also that $\norm{a} = \norm{a^*}$ for all $a$.
We will frequently use that
$a^2 = -\norm{a}^2$ 
if $a$ is an imaginary element of $A_n$.

\begin{lemma}
\label{lem:C-bi-conj2}
Let $a$ belong to $\C_n^\perp$, and let $\alpha$ and $\beta$
belong to $\C_n$.  Then $\alpha a \cdot \beta a = -\norm{a}^2 \alpha \beta^*$.
\end{lemma}

\begin{proof}
Follow the same general strategy as in the proof of Lemma \ref{lem:C-bi-conj}.
However, instead of using that $ab$ belongs to $\C_n^\perp$,
use that $a^2 = -\norm{a}^2$ is real.
\end{proof}

One consequence of Lemma \ref{lem:C-bi-conj2} is that 
$\norm{\alpha a} = \norm{\alpha} \norm{a}$ if
$\alpha$ belongs to $\C_n$ and $a$ belongs to $\C_n^\perp$.
This follows from the computation 
$\alpha a \cdot \alpha a = -\norm{a}^2 \alpha \alpha^*$.

%%%%%%%%%%%%%%%%%%%%%%%%%%%%%%%%%%%%%%%%%%%%%%%%%%%%%%%%%%%

\subsection{The subalgebra $\HH_n$}

\begin{defn}
Let \mdfn{$\HH_n$} be the $\R$-linear span of the elements
$1$, $i_n$, $i_{n-1}$, and $i_{n-1} i_n$ of $A_n$.
\end{defn}

The notation reminds us that $\HH_n$ is a subalgebra isomorphic to
the quaternions.  
Many of the results that follow
refer to $\HH_n$ and its orthogonal complement $\HH_n^\perp$.

In terms of the product
$A_n = A_{n-1} \times A_{n-1}$, $\HH_n$ is the $\R$-linear span
of $(1,0)$, $(0,1)$, $(i_{n-1}, 0)$, and $(0,i_{n-1})$.
By inspection, $\HH_n$ is a $\C_n$-linear subspace of $A_n$.
It is also equal to $\C_{n-1} \times \C_{n-1}$.
Also, $\HH_n^\perp$ and $\C_{n-1}^\perp \times \C_{n-1}^\perp$
are equal as subspaces of $A_n$.

%%%%%%%%%%%%%%%%%%%%%%%%%%%%%%%%%%%%%%%%%%%%%%%%%%%%%%%%%%%%%%%%%%%%

\subsection{Zero-divisors and annihilators}

A \mdfn{zero-divisor} is a non-zero element $a$ of $A_n$ such that
there exists another non-zero element $b$ in $A_n$ with $ab = 0$.
The \mdfn{annihilator $\Ann(a)$}
of $a$ is the set of all elements $b$ such that $ab = 0$.
In other words, $\Ann(a)$ is the kernel of left multiplication by $a$.

\begin{lemma}[Corollary 1.9 of M1 and Lemma 9.5 of DDD]
\label{lem:zd-C-perp}
If $a$ is a zero-divisor in $A_n$, then $a$ belongs to $\C_n^\perp$.
\end{lemma}

\begin{thm}[Theorem 9.8 and Proposition 9.10 of DDD]
\label{thm:4dim}
The dimension of any annihilator in $A_n$ is a multiple of 4
and is at most $2^n - 4n + 4$.
\end{thm}

See also \cite[Corollary 1.17]{M1} for another proof of the first claim.

\begin{lemma}
\label{lem:ann-im}
Let $a$ belong to $\C_n^\perp$.
For any $b$ in $A_n$, the product $ab$ is orthogonal to $\Ann(a)$.
\end{lemma}

\begin{proof}
Let $c$ belong to $\Ann(a)$.  Use Lemma \ref{lem:C-conj-linear} 
to deduce that
$\langle ab, c \rangle_{\C} = -\langle b, ac \rangle^*_{\C}$.
This equals zero because $ac = 0$.
\end{proof}

Let \mdfn{$\Im(a)$} be the image of left multiplication by $a$.
Lemma \ref{lem:ann-im} implies that $\Im(a)$ is the orthogonal complement
of $\Ann(a)$ in $A_n$.

%%%%%%%%%%%%%%%%%%%%%%%%%%%%%%%%%%%%%%%%%%%%%%%%%%%%%%%%%%%%%%%%

\subsection{Projections}
\label{subsctn:proj}

We still need a few technical definitions and results.  We provide
complete proofs for the following results because their proofs do not
already appear elsewhere.

\begin{defn}
\label{defn:proj}
For any $a$ in $A_n$, let \mdfn{$\pi_{\C}(a)$} be the orthogonal
projection of $a$ onto $\C_n$, and let
\mdfn{$\pi_{\C}^\perp (a)$} be the 
orthogonal projection of $a$ onto $\C_n^\perp$.
\end{defn}

By definition,
$\pi_{\C}(ab^*)$ equals $\langle a, b \rangle_{\C}$
for any $a$ and $b$.

\begin{lemma}
\label{lem:proj-multiply}
Let $a$ and $b$ belong to $A_n$. Let $b = b' + b''$,
where $b'$ is the $\C$-orthogonal projection of $b$ onto the $\C$-linear
span of $a$ and 
where $b''$ is the $\C$-orthogonal projection of $b$ onto 
the $\C$-orthogonal complement of $a$.
Then $\pi_{\C}(ab) = ab'$, and
$\pi_{\C}^\perp(ab) = ab''$.
Similarly,
$\pi_{\C}(ba) = b'a$, and
$\pi_{\C}^\perp(ba) = b''a$.
\end{lemma}

\begin{proof}
Note that $ab = ab' + ab''$.  The first term belongs
to $\C_n$ by Lemma \ref{lem:C-bi-conj2}, 
while the second term belongs to $\C_n^\perp$
because $a$ and $b''$ are $\C$-orthogonal.
Similarly, $ba = b'a + b''a$, where
$b'a$ belongs to $\C_n$ and $b''a$ belongs to $\C_n^\perp$.
\end{proof}

\begin{cor}
\label{cor:C-proj}
For any $a$ and $b$ in $\C_n^\perp$, 
$\pi_{\C}(ab) = \pi_{\C}(ba)^*$ and
$\pi_{\C}^\perp (ab) = -\pi_{\C}^\perp(ba)$.
\end{cor}

\begin{proof}
Write $b = b' + b''$,
where $b'$ is the $\C$-orthogonal projection of $b$ onto $a$ and 
where $b''$ is the $\C$-orthogonal projection of $b$ onto 
the $\C$-orthogonal complement of $a$.
By Lemma \ref{lem:proj-multiply}, $\pi_{\C}(ab) = ab'$ and
$\pi_{\C}(ba) = b'a$.  It follows from Lemma \ref{lem:C-bi-conj2} that
$(ab')^* = b'a$.  This finishes the first claim.

For the second claim,
Lemma \ref{lem:proj-multiply} implies that $\pi_{\C}^\perp(ab) = ab''$ and
$\pi_{\C}^\perp(ba) = b''a$.  
Because $a$ and $b''$ are imaginary and orthogonal,
$ab'' = -b''a$.
\end{proof}

\begin{cor}
\label{cor:proj-multiply}
Let $a$ belong to $A_n$, and let $\alpha$ belong to $\C_n$.
Then $\pi_{\C}(\alpha a) = \alpha \pi_{\C}(a) = \pi_{\C}(a \alpha)$.
\end{cor}

\begin{proof}
This is an immediate consequence of Lemma \ref{lem:proj-multiply}
and the fact that $\C_n$ is commutative.
\end{proof}

One way to interpret Corollary \ref{cor:proj-multiply} is that
$\pi_{\C}$ is a $\C$-linear map.

\begin{cor}
\label{cor:C-multiply}
Let $a$ and $b$ belong to $A_n$.
Then $ab$ belongs to $\C_n$ if and only if 
$b$ belongs to the $\C$-linear span of $a$ and $\Ann(a)$.
\end{cor}

\begin{proof}
In the notation of Lemma \ref{lem:proj-multiply}, observe
that $ab$ belongs to $\C_n$ if and only if $ab''$ is zero.
\end{proof}

%%%%%%%%%%%%%%%%%%%%%%%%%%%%%%%%%%%%%%%%%%%%%%%%%%%%%%%%%%%

\section{Notation}
\label{sctn:bracket}

\begin{defn}
\label{defn:bracket}
For any $a$ and $b$ in $\C_n^\perp$, let \mdfn{$\{a,b\}$} be the element
\[
\frac{1}{\sqrt{2}}\big( a+b, i_n(-a+b) \big)
\]
of $A_{n+1}$.
\end{defn}

Whenever we write an expression of the form $\{a, b\}$, the reader
should automatically
assume that $a$ and $b$ belong to $\C_n^\perp$; nevertheless,
we have tried to be explicit with this assumption.
The reason for the factors $\frac{1}{\sqrt{2}}$ will show up 
in Lemma \ref{lem:convert-norm} and
Lemma \ref{lem:bracket-inner}, where we study 
the metric properties of the notation $\{ a, b\}$.

\begin{lemma}
\label{lem:convert}
Let $(x,y)$ belong to $\HH_{n+1}^\perp$, i.e., let $x$ and $y$ belong to 
$\C_n^\perp$.  Then
\[
(x,y) = \frac{1}{\sqrt{2}} \big\{x+ i_n y, x-i_n y \big\}.
\]
The subspace $\HH_{n+1}^\perp$ of $A_{n+1}$ is equal to the subspace of all
elements of the form $\{a, b\}$ with $a$ and $b$ in $\C_n^\perp$.
\end{lemma}

\begin{proof}
For the first claim, check the definition.  
This immediately implies that
every element of $\HH_{n+1}^\perp$
can be written in the form $\{a, b\}$ for some $a$ and $b$ in $\C_n^\perp$.

On the other hand,
let $a$ and $b$ belong to $\C_n^\perp$.
Then $a+b$ and $i_n(-a+b)$ also belong to $\C_n^\perp$,
so $\big(a+b, i_n(-a+b) \big)$ belongs to $\HH_{n+1}^\perp$.
\end{proof}

Recall that left multiplication makes $A_{n+1}$ into a
$\C_{n+1}$-vector space.  We now describe multiplication by elements
$\C_{n+1}$ with respect to the notation 
$\{a, b\}$.  

\begin{defn}
\label{defn:tilde}
If $\alpha$ belongs to $\C_n$, then \mdfn{$\tilde{\alpha}$}
is the image of $\alpha$ under the 
$\R$-linear map $\C_n \map \C_{n+1}$
that takes 1 to 1 and $i_n$ to $i_{n+1}$.  
\end{defn}

\begin{lemma}
\label{lem:bracket-C-action}
Let $a$ and $b$ belong to $\C_n^\perp$, and let $\alpha$ belong to
$\C_n$.  Then
\[
\tilde{\alpha} \{a, b\} = \{ \alpha^* a, \alpha b\}.
\]
\end{lemma}

\begin{proof}
Compute directly that $i_{n+1} \{a, 0 \} = \{-i_n a, 0\}$ and
$i_{n+1} \{0, b \} = \{0, i_n b\}$.
\end{proof}

\begin{lemma}
\label{lem:convert-norm}
For any $a$ and $b$ in $\C_n^\perp$,
\[
\norm{\{a, b\}}^2 = \norm{a}^2 + \norm{b}^2.
\]
\end{lemma}

\begin{proof}
According to Definition \ref{defn:bracket}, $\norm{\{a,b\}}^2$
equals 
\[
\frac{1}{2} \big( \norm{a+b}^2 + \norm{i_n(-a+b)}^2 \big).
\]
As a consequence of Lemma \ref{lem:C-bi-conj2}, this expression equals
\[
\frac{1}{2} \big( \norm{a+b}^2 + \norm{-a+b}^2 \big),
\]
which simplifies to $\norm{a}^2 + \norm{b}^2$ by the parallelogram law.
\end{proof}

The absence of scalars in the above formula
is the primary reason that
the scalar  $\frac{1}{\sqrt{2}}$ appear in Definition \ref{defn:bracket}.

%%%%%%%%%%%%%%%%%%%%%%%%%%%%%%%%%%%%%%%%%%%%%%%%%%%%%%%%%%%%%%%%%%%

\section{Multiplication Formulas}
\label{sctn:mult}

This section is the technical heart of the paper.  We establish formulas
for multiplication with respect to the notation of Section \ref{sctn:bracket}.
The rest of the paper consists of many applications of these formulas.

\begin{prop}
\label{prop:bracket-multiply}
Let $a$, $b$, $x$, and $y$ belong to $\C_n^\perp$, and suppose that
$a$ and $b$ are both $\C$-orthogonal to $x$ and $y$.
Then
\[
\{a, b\} \{x, y \} = \sqrt{2} \{ ax, by \}.
\]
\end{prop}

\begin{proof}
We begin by computing that $\{a, 0\}\{x,0\}$ equals
\[
\frac{1}{2}\Big( ax + i_n x \cdot i_n a, -i_n x \cdot a + i_n a \cdot x \Big).
\]
Apply Lemma \ref{lem:C-bi-conj} to simplify this expression to
\[
\frac{1}{2}\Big( ax - xa, i_n \cdot xa - i_n \cdot ax \Big).
\]
Note that $ax = -xa$ because $a$ and $x$ are imaginary and orthogonal,
so this expression further simplifies to
\[
( ax , - i_n \cdot ax ),
\]
which equals $\sqrt{2}\{ax, 0\}$.
A similar calculation shows that
\[
\{0, b\}\{0,y\} = \sqrt{2} \{0, by\}.
\]

Next we compute that $\{a, 0\}\{0,y\}$ equals
\[
\frac{1}{2} \Big( ay - i_n y \cdot i_n a, i_n y \cdot a + i_n a \cdot y \Big).
\]
Again use Lemma \ref{lem:C-bi-conj} to simplify to
\[
\frac{1}{2} \Big( ay + ya, -i_n \cdot ya - i_n \cdot ay \Big),
\]
but this equals zero because $ay = -ya$.

A similar calculation shows that
$\{0, b\}\{x,0\} = 0$.
\end{proof}

\begin{remark}
Proposition \ref{prop:bracket-multiply} already
gives a sense of how easy it is
to express certain zero-divisors using the notation $\{a,b\}$.
For example, the product $\{a,0\}\{0,y\}$ 
is always zero as long as $a$ and $y$ are $\C$-orthogonal elements
of $\C_n^\perp$.
\end{remark}

Because of the orthogonality hypotheses on $a$, $b$, $x$, and $y$,
Proposition \ref{prop:bracket-multiply} does not quite describe
how to multiply arbitrary elements of $\HH_{n+1}^\perp$.  Therefore,
we need more multiplication formulas to handle various special cases.

\begin{lemma}
\label{lem:bracket-multiply-parallel}
Let $a$ belong to $\C_n^\perp$.  Then
\[
\{0, a \} \{a, 0\} = - \{a, 0\} \{0,a\} = \norm{a}^2(0,i_n).
\]
\end{lemma}

\begin{proof}
Compute that $\{0, a\}\{a,0\}$ equals
\[
\textstyle \frac{1}{2} \Big( a^2 - i_n a \cdot i_n a, -2i_n a \cdot a \Big).
\]
Lemma \ref{lem:C-bi-conj2} implies
that the first coordinate is zero and that
the second coordinate is $\norm{a}^2 i_n$.

Finally, observe that $\{0,a\}$ and $\{a,0\}$ are orthogonal
and imaginary; therefore they anti-commute.
\end{proof}

We write \mdfn{$\tilde{\pi}_{\C}$} for the composition of the 
projection $A_n \map \C_n$ with the map $\C_n \map \C_{n+1}$ described
in Definition \ref{defn:tilde}.

\begin{cor}
\label{cor:bracket-multiply-parallel}
Let $a$ and $b$ be $\C_n$-linearly dependent elements of $\C_n^\perp$.
Then
\begin{enumerate}
\item
$\{ a, 0 \} \{ b, 0 \}  =  \tilde{\pi}_{\C} (ab)^*$.
\item
$\{ 0, a \} \{ 0, b \} = \tilde{\pi}_{\C} (ab)$.
\item
$\{ a, 0\} \{ 0, b \} = \tilde{\pi}_{\C} (ab) \cdot (0,i_n)$.
\item
$\{ 0, a \} \{ b, 0 \} = -\tilde{\pi}_{\C} (ab)^* \cdot (0,i_n)$.
\end{enumerate}
\end{cor}

\begin{proof}
Since $b$ belongs to the $\C_n$-linear span of $a$, we may write
$b = \alpha a$ for some $\alpha$ in $\C_n$.  
Lemma \ref{lem:C-bi-conj2} implies that $ab$ equals $-\norm{a}^2 \alpha^*$,
so $\tilde{\pi}_{\C}(ab)$ equals $-\norm{a}^2 \tilde{\alpha}^*$.

On the other hand, $\{a,0\}\{\alpha a, 0\}$ equals
$\{a,0\} \cdot \tilde{\alpha}^* \{a,0\}$ by Lemma \ref{lem:bracket-C-action},
which also equals $-\norm{\{a,0\}}^2 \tilde{\alpha}$
by Lemma \ref{lem:C-bi-conj2}.  Finally,
this equals $-\norm{a}^2 \tilde{\alpha}$ by 
Lemma \ref{lem:convert-norm}.  This establishes formula (1).
The calculation for formula (2) is similar.

Next, $\{a,0\}\{0,\alpha a\}$ equals
$\{a,0\} \cdot \tilde{\alpha} \{0,a\}$ by Lemma \ref{lem:bracket-C-action},
which also equals $\tilde{\alpha}^* \cdot \{a,0\}\{0,a\}$
by Lemma \ref{lem:C-bi-conj}.  Finally,
this equals $-\norm{a}^2 \tilde{\alpha}^*(0,i_n)$ by 
Lemma \ref{lem:bracket-multiply-parallel},
establishing formula (3).
The calculation for formula (4) is similar.
\end{proof}

We are now ready to give an explicit formula for multiplication
of arbitrary elements of $\HH_{n+1}^\perp$.

\begin{thm}
\label{thm:bracket-multiply}
Let $a$, $b$, $x$, and $y$ belong to $\C_n^\perp$.  Then
$\{ a, b \} \{x, y \}$ equals
\[
\sqrt{2} \big\{ \pi_{\C}^\perp (ax), \pi_{\C}^\perp (by) \big\} +
\tilde{\pi}_{\C} \big( xa + by \big) + 
\tilde{\pi}_{\C} \big( ay - xb \big) (0,i_n).
\]
\end{thm}

\begin{proof}
We begin by computing $\{a, 0\} \{x, 0\}$.
Write $x = x' + x''$, where $x'$ belongs to the $\C$-linear span of $a$ and
$x''$ is $\C$-orthogonal to $a$.  Then 
\[
\{a, 0\} \{x, 0 \} =
\{a, 0\} \{x', 0\} + \{a,0\} \{x'', 0\}.
\]
The first term equals $\tilde{\pi}_{\C}(ax')^*$ by 
Corollary \ref{cor:bracket-multiply-parallel},
which in turn equals $\tilde{\pi}_{\C}(x'a)$ by 
Corollary \ref{cor:C-proj}.
This is the same as $\tilde{\pi}_{\C} (xa)$ by 
Lemma \ref{lem:proj-multiply}.
The second term equals $\sqrt{2}\{ax'', 0\}$
by Proposition \ref{prop:bracket-multiply}, which equals
$\sqrt{2} \{ \pi_{\C}^\perp(ax), 0 \}$ by Lemma \ref{lem:proj-multiply}.
The computation for $\{0, b\} \{0, y\}$ is similar.

Now consider the product $\{a, 0\} \{0, y\}$.
Write $y = y' + y''$, where $y'$ belongs to the $\C$-linear span of $a$ and
$y''$ is $\C$-orthogonal to $a$.  Then 
\[
\{a, 0\} \{0, y \} =
\{a, 0\} \{0, y' \} + \{a,0\} \{0,y''\}.
\]
The first term equals $\tilde{\pi}_{\C}(ay') \cdot (0,i_n)$ by 
Corollary \ref{cor:bracket-multiply-parallel},
which is the same as $\tilde{\pi}_{\C} (ay) \cdot (0,i_n)$ by 
Lemma \ref{lem:proj-multiply}.
The second term equals zero
by Proposition \ref{prop:bracket-multiply}.
The computation for $\{0, b\} \{x, 0\}$ is similar.
\end{proof}

\begin{remark}
\label{rem:bracket-multiply}
The three terms in the formula of Theorem \ref{thm:bracket-multiply}
are orthogonal.  The first term belongs to $\HH_{n+1}^\perp$;
the second term belongs to $\C_{n+1}$; and the third term
belongs to $\HH_{n+1} \cap \C_{n+1}^\perp$, which is also the
$\C$-linear span of $(0,i_n)$ or the $\R$-linear span 
of $(i_n,0)$ and $(0,i_n)$.
\end{remark}

Theorem \ref{thm:bracket-multiply} shows how to compute the product
of two elements of $\HH_{n+1}^\perp$.  
On the other hand, 
it is easy to multiply elements of $\HH_{n+1}$; this is just ordinary
quaternionic arithmetic.  
In order to have a complete description of multiplication on $A_{n+1}$,
we need to explain how to multiply elements of $\HH_{n+1}$ 
with elements of $\HH_{n+1}^\perp$.

Lemma \ref{lem:bracket-C-action} shows how to compute the product
of an element of $\HH_{n+1}^\perp$ and an element of $\C_{n+1}$.
It remains only to compute the product 
of an element of $\HH_{n+1}^\perp$ and an element of 
$\HH_{n+1} \cap \C_{n+1}^\perp$, i.e., 
the $\C$-linear span of $(0,i_n)$.
The following lemma makes this computation.

\begin{lemma}
\label{lem:last-multiply}
Let $a$ and $b$ belong to $\C_n^\perp$.  Then
\[
(0,i_n)\{a,b\} = 
-\{a,b\} (0,i_n) = 
\{b, -a\}.
\]
\end{lemma}

\begin{proof}
Compute directly that $(0,i_n)\{a,0\} = \{0,-a\}$ and that
$(0,i_n)\{0,b\} = \{b,0\}$.
Also, use that orthogonal imaginary elements anti-commute.
\end{proof}

\subsection{Inner product computations}

\begin{lemma}
\label{lem:bracket-inner}
Let $a$, $b$, $x$, and $y$ belong to $\C_n^\perp$.  Then
\[
\big\langle \{a,b\},\{x,y\} \big\rangle_{\C} = 
 \langle a,x\rangle^*_{\C} + 
                       \langle b,y\rangle_{\C}.
\]
\end{lemma}

\begin{proof}
We need to compute the projection of the product
$-\{a, b\} \{x,y\}$ onto $\C_{n+1}$.  
Theorem \ref{thm:bracket-multiply} immediately shows that this 
projection equals $-\tilde{\pi}_{\C} (xa + by)$,
which is equal to 
$\langle x, a \rangle_{\C} + \langle b,y\rangle_{\C}$.
Finally, recall that $\langle x, a \rangle_{\C} = \langle a, x \rangle^*_{\C}$.
\end{proof}

\begin{cor}
\label{cor:bracket-inner}
Let $a$, $b$, $x$, and $y$ belong to $\C_n^\perp$.  Then
\[
\big\langle \{a,b\},\{x,y\} \big\rangle_{\R} = 
\langle a,x\rangle_{\R} +\langle b,y\rangle_{\R}.
\]
\end{cor}

\begin{proof}
Use Lemma \ref{lem:bracket-inner}, 
recalling that the real inner product equals
the real part of the Hermitian inner product.  
\end{proof}

\subsection{Subalgebras}
\label{subsctn:eigenpair}

Suppose that $a$ and $b$ are $\C$-orthogonal elements of $\C_n^\perp$
that both have norm $1$.  Suppose also that $a$ and $b$
satisfy the equations 
$a \cdot ab = -\lambda b$ and $b \cdot ba = -\lambda a$ for some
non-zero real number $\lambda$.  These equations 
guarantee that $a$ and $b$ generate a $4$-dimensional
subalgebra of $A_n$; the subalgebra is isomorphic to $\HH$ when $\lambda = 1$.
This remark concerns the possible values for $\lambda$, and therefore 
addresses 
the classification problem for 4-dimensional subalgebras of Cayley-Dickson
algebras.
See \cite[Section~7]{CD} 
for detailed information on 4-dimensional subalgebras
of $A_4$.  In particular, in $A_4$, the only possible values for 
$\lambda$ are $1$ and $2$ \cite[Theorem~7.1]{CD}.

Given $a$ and $b$ as in the previous paragraph,
compute that 
\[
\frac{1}{\sqrt{2}}\{a,b\} \cdot \frac{1}{\sqrt{2}} \{b, -a\} = 
\frac{1}{\sqrt{2}}\{ab, -ba\} + (0,i_n)
\]
using Theorem \ref{thm:bracket-multiply}.
Next, compute that 
\[
\frac{1}{\sqrt{2}}\{a,b\} 
\left( \frac{1}{\sqrt{2}}\{ab,-ba\} + (0,i_n) \right) =
-\frac{\lambda + 1}{\sqrt{2}} \{b,-a\} 
\]
using Proposition \ref{prop:bracket-multiply} and 
Lemma \ref{lem:last-multiply}.
This uses that $a$ and $ab$ are $\C$-orthogonal by 
Lemma \ref{lem:ortho1}
and also the equations involving $a$, $b$, and $\lambda$.
A similar calculation can be performed with the roles of
$\frac{1}{\sqrt{2}}\{a,b\}$ and $\frac{1}{\sqrt{2}}\{b,-a\}$ switched.

We have shown that $\frac{1}{\sqrt{2}} \{a,b\}$ and
$\frac{1}{\sqrt{2}} \{b,-a\}$ satisfy the same equations as
$a$ and $b$ do, except that $\lambda$ is replaced by $\lambda +1$.
Using the argument of \cite[Theorem~7.1]{CD} (which can be applied
even when $n > 4$), it follows that for every positive integer $r$
and every sufficiently large $n$ (depending on $r$),
there is a subalgebra of $A_n$ that is isomorphic to the non-associative
algebra with $\R$-basis $\{1, x, y, z\}$ subject to the multiplication
rules
\[
x^2 = y^2 = z^2 = -1, \quad xy = -yx = z\sqrt{r}, \quad
yz = -zy = x\sqrt{r}, \quad zx = -xz = y\sqrt{r}.
\]
This algebra is isomorphic to $\HH$ when $r = 1$.

Another consequence of our multiplication formulas is the
following observation about sets of mutually annihilating elements.

\begin{lemma}
\label{lem:Zm}
Let $n \geq 3$.
If $\C_n^\perp$ contains two sets
$\{ x_1, \ldots, x_{2^{n-3}} \}$ and 
$\{ y_1, \ldots, y_{2^{n-3}} \}$ of size $2^{n-3}$ such that
$x_i x_j = 0 = y_i y_j$ for all $i \neq j$
and each $x_i$ is $\C$-orthogonal to each $y_j$, then the product
$\{x_i,0\}\{x_j,0\}$ is zero when $i \neq j$, and 
$\{x_i,0\}\{0,y_j\}$ is zero for all $i$ and $j$.
\end{lemma}

\begin{proof}
Compute with Proposition \ref{prop:bracket-multiply}.
\end{proof}

\begin{cor}
\label{cor:Zm}
The space $\C_n^\perp$ contains $2^{n-3}$ elements 
such that the product of any two distinct elements is zero.
\end{cor}

\begin{proof}
We will actually prove a stronger result that
$\C_n^\perp$ contains two sets
$\{ x_1, \ldots, x_{2^{n-3}} \}$ and 
$\{ y_1, \ldots, y_{2^{n-3}} \}$ of size $2^{n-3}$ such that
$x_i x_j = 0 = y_i y_j$ for all $i \neq j$
and each $x_i$ is $\C$-orthogonal to each $y_j$.

The proof is by induction on $n$, using Lemma \ref{lem:Zm}.
The base case $n=3$ is trivial; it just calls for the
existence of two orthogonal elements of the six-dimensional subspace 
$\C_3^\perp$ of $A_3$.

Now suppose for induction that the sets $\{x_1, \ldots, x_{2^{n-3}} \}$
and $\{y_1, \ldots, y_{2^{n-3}} \}$ exist in $A_n$.
Consider the subset of $A_{n+1}$ consisting of all elements of the 
form $\{x_i,0\}$ or $\{0,y_j\}$.  There are $2^{n-2}$ such elements,
and Lemma \ref{lem:Zm} implies that the product of any two distinct
such elements is zero.

Also consider the subset of $A_{n+1}$ consisting of all elements of the 
form $\{y_j,0\}$ or $\{0,x_i\}$.  Again, there are $2^{n-2}$ such elements,
and the product of any two distinct such elements is zero.

Finally, by Lemma \ref{lem:bracket-inner} 
and the induction assumption, the elements
described in the previous two paragraphs are $\C$-orthogonal.
\end{proof}

Corollary \ref{cor:Zm}
is also relevant to subalgebras of Cayley-Dickson algebras.
The $\R$-linear span of $1$ together with a set of mutually
annihilating elements is a subalgebra of $A_n$.  These subalgebras are
highly degenerate in the sense that $xy = 0$ for any pair of orthogonal
imaginary elements.  Corollary \ref{cor:Zm} implies that $A_n$ contains
such a subalgebra of dimension $1 + 2^{n-3}$.  In fact, we have shown
that $A_n$ contains two such subalgebras whose imaginary
parts are $\C$-orthogonal.

\begin{ques}
Does $A_n$ contain a degenerate subalgebra of dimension larger
than $1 + 2^{n-3}$?
\end{ques}

%%%%%%%%%%%%%%%%%%%%%%%%%%%%%%%%%%%%%%%%%%%%%%%%%%%%%%%%%%%%%%%%%%

\section{Annihilation in $\HH_{n+1}^\perp$}
\label{sctn:ann-H-perp}

In this section, we apply the multiplication formulas of 
Section \ref{sctn:mult} to consider zero-divisors in $A_{n+1}$.

\begin{prop}
\label{prop:bracket-zd}
Let $a$, $b$, $x$, and $y$ belong to $\C_{n}^\perp$.  
Then $\{a,b\}\{x,y\}=0$ if and only if 
\begin{enumerate}[(i)]
\item $\pi_{\C}^\perp(ax)=0$,
\item $\pi_\C^\perp(by)=0$,
\item $xa+by=0$, and
\item $\pi_\C(ay-xb)=0$.
\end{enumerate}
\end{prop}

\begin{proof}
Parts (i), (ii), and (iv) 
are immediate from Theorem \ref{thm:bracket-multiply}.
It follows from (i) and (ii) that $\pi_\C(xa+by)=xa+by$.
Therefore, part (iii) also follows from 
Theorem \ref{thm:bracket-multiply}.
\end{proof}

The conditions of Proposition \ref{prop:bracket-zd} are redundant.
For example, condition (i) follows from conditions (ii) and (iii).
However, it is more convenient to formulate the proposition symmetrically.

\begin{prop}
\label{prop:ann-intersect-H-perp}
Let $n \geq 3$.
Let $a$ and $b$ be non-zero elements of $\C_n^\perp$.
Then
$\HH_{n+1}^\perp \cap \Ann\{a,b\}$ is equal to the space of all
$\{ \alpha a + x, \beta b + y\}$ such that:
\begin{enumerate}
\item 
$x$ belongs to $\Ann(a)$, and
$y$ belongs to $\Ann(b)$;
\item
$\alpha$ and $\beta$  belong to $\C_{n}$;
\item
$\norm{a}^2 \alpha + \norm{b}^2 \beta^* = 0$; and
\item
$(\beta^* - \alpha)\pi_{\C}(ab) + \pi_{\C}(ay-xb) = 0$.
\end{enumerate} 
\end{prop}

\begin{proof}
We want to solve the equation
$\{a,b\}\{z,w\} = \{0,0\}$ under the assumption that
$z$ and $w$ belong to $\C_n^\perp$ (see Lemma \ref{lem:convert}).
Using Proposition~\ref{prop:bracket-zd}, this is equivalent to
solving the four equations
\begin{myequation}
\label{eq:aiHp1}
\pi_{\C}^\perp(az) & = 0 \\
  \addtocounter{subsection}{1}
\label{eq:aiHp2}
\pi_{\C}^\perp(bw) & = 0 \\
  \addtocounter{subsection}{1}
\label{eq:aiHp3}
za+bw & = 0 \\
  \addtocounter{subsection}{1}
\label{eq:aiHp4}
\pi_{\C}(aw-zb) & = 0.
\end{myequation}
By Corollary~\ref{cor:C-multiply}, Equations (\ref{eq:aiHp1}) and
(\ref{eq:aiHp2}) are the same
as requiring that $z$ belongs to 
the $\C$-linear span of $a$ and $\Ann(a)$
and that $w$ belongs to 
the $\C$-linear span of $b$ and $\Ann(b)$.
Therefore, we may write
$z = \alpha a + x$ and $w = \beta b + y$ for some
$\alpha$ and $\beta$ in $\C_n$, some $x$ in $\Ann a$, and some
$y$ in $\Ann b$.  We also know that $x$ and $y$ belong to $\C_n^\perp$
by Lemma \ref{lem:zd-C-perp}; this is where we use that $a$ and $b$
are non-zero.

Substitute the expressions for $z$ and $w$ in Equations
(\ref{eq:aiHp3}) and (\ref{eq:aiHp4}) to obtain 
\begin{myequation}
\label{eq:aiHp5}
(\alpha a + x)a+ b(\beta b + y) & = 0 \\
  \addtocounter{subsection}{1}
\label{eq:aiHp6}
\pi_{\C}\Big(a(\beta b + y) - (\alpha a + x)b \Big) & = 0.
\end{myequation}
Equation (\ref{eq:aiHp5}) 
simplifies to $-\norm{a}^2\alpha - \norm{b}^2 \beta^* = 0$
by Lemma \ref{lem:C-bi-conj2} and the fact that $xa = by = 0$.
This is condition (3) of the proposition.

Equation (\ref{eq:aiHp6}) can be rewritten as 
\begin{myequation}
\label{eq:aiHp7}
\pi_{\C}( \beta^* \cdot ab - ab \cdot \alpha ) +
\pi_{\C} ( ay - xb ) = 0
\end{myequation}
by Lemma \ref{lem:C-conj-linear}.
Apply Corollary \ref{cor:proj-multiply} to the second part of the
first term of Equation (\ref{eq:aiHp7}) to obtain the equation
$(\beta^* - \alpha )\pi_{\C}(ab) + \pi_{\C}( ay - xb ) = 0$.
This is condition (4) of the proposition.
\end{proof}

\begin{thm}
\label{thm:ann-bracket-bound}
Let $n \geq 3$, and 
let $a$ and $b$ be non-zero elements of $\C_{n}^\perp$.  Then 
$\dim \Ann\{a,b\}$ equals  $\dim \Ann a + \dim \Ann b$ or
$\dim \Ann a + \dim \Ann b +4$.
\end{thm}

\begin{proof}
First we will use Proposition \ref{prop:ann-intersect-H-perp}
to analyze $\HH_{n+1}^\perp \cap \Ann\{a,b\}$.
Let $V$ be the space of elements $\{ \alpha a + x, \beta b + y\}$
such that $\alpha$ and $\beta$
belong to $\C_n$, $x$ belongs to $\Ann a$, and $y$ belongs to $\Ann b$.
The dimension of $V$ is equal to 
$\dim \Ann a + \dim \Ann b + 4$.
Recall from Lemma \ref{lem:bracket-C-action} that
for $\gamma$ in $\C_n$,
\[
\tilde{\gamma} \{\alpha a + x, \beta b + y \} = 
\{\gamma^* \alpha a + \gamma^* x, \gamma \beta b + \gamma y \}.
\]
This shows that $V$ is a $\C_n$-vector space, and
Condition (3) of Proposition \ref{prop:ann-intersect-H-perp}
is a non-degenerate conjugate-linear 
equation in the variables $\alpha$ and $\beta$.
Hence there is a subspace
of $V$ of dimension $\dim \Ann a + \dim \Ann b + 2$ that satisfies
condition (3).  

Condition (4) of Proposition \ref{prop:ann-intersect-H-perp}
is a conjugate-linear equation in the variables $\alpha$, $\beta$,
$x$, and $y$, which 
may or may not be non-degenerate
and independent of condition (3).  This establishes that
\[ 
\dim \Ann a+\dim \Ann b \leq \dim (\HH_{n+1}^\perp \cap \Ann\{a,b\})
\leq \dim \Ann a +\dim \Ann b +2. 
\]

Lemma \ref{lem:zd-C-perp} implies that $\Ann\{a,b\}$ is contained
in $\C_{n+1}^\perp$.
Note that $\HH_{n+1}^\perp$ is a codimension 2 subspace of $\C_{n+1}^\perp$.
Therefore, 
the codimension of $\HH_{n+1}^\perp \cap \Ann\{a,b\}$ in $\Ann\{a,b\}$
is at most 2.  This establishes the inequality
\[
\dim \Ann a+\dim \Ann b \leq \dim \Ann\{a,b\}
\leq \dim \Ann a +\dim \Ann b +4. 
\]
The desired result follows from Theorem \ref{thm:4dim}, which tells
us that the dimension of any annihilator is a multiple of 4.
\end{proof}

Theorem \ref{thm:ann-bracket-bound} gives two options for the
dimension of $\Ann\{a,b\}$; Section \ref{sctn:D-locus} below contains
conditions on $a$ and $b$ that distinguish between these two cases.

One might also be concerned that Theorem \ref{thm:ann-bracket-bound}
applies only  to elements $\{a,b\}$ in which
both $a$ and $b$ are non-zero
because it relies on Proposition \ref{prop:ann-intersect-H-perp}.
For completeness, we also review from \cite{DDD}
the simpler situation of elements of the form $\{a,0\}$ and $\{0,a\}$.
The following proposition 
can be proved with the formulas of Section \ref{sctn:mult}.

\begin{prop}[Theorem 10.2, DDD]
\label{prop:ann-bracket-special}
Let $n \geq 4$, and let $a$ belong to $\C_{n-1}^\perp$.  Then 
the element $\{a,0\}$ of $A_n$ is a zero-divisor whose annihilator
$\Ann\{a,0\}$ equals the space of all elements $\{x,y\}$
where $x$ belongs to $\Ann(a)$ and 
$y$ is $\C$-orthogonal to $1$ and $a$.
Similarly,
the element $\{0,a\}$ of $A_n$ is a zero-divisor whose annihilator
$\Ann\{0,a\}$ equals the space of all elements $\{x,y\}$
where 
$y$ belongs to $\Ann(a)$ and
$x$ is $\C$-orthogonal to $1$ and $a$.
In either case, the dimension of the annihilator is
$\dim \Ann(a)+2^{n-1}-4$.
\end{prop}

In fact, \cite[Theorem 10.2]{DDD} was a major inspiration for the notation
$\{a,b\}$.

%%%%%%%%%%%%%%%%%%%%%%%%%%%%%%%%%%%%%%%%%%%%%%%%%%%%%%%%%%%%%%%%

\section{The $D$-locus}
\label{sctn:D-locus}

In Section \ref{sctn:ann-H-perp},
we started to consider $\Ann\{a,b\}$ when $a$
and $b$ are arbitrary elements in $\C_n^\perp$, i.e., when 
$\{a,b\}$ is an arbitrary element of $\HH_{n+1}^\perp$.
Theorem \ref{thm:ann-bracket-bound} told us that
except for some simple well-understood cases covered in 
Proposition \ref{prop:ann-bracket-special}, the
dimension of $\Ann\{a,b\}$ is either $\dim \Ann a+\dim \Ann b$ or 
$\dim \Ann a+\dim \Ann b+4$.  The goal of this section is to
distinguish between these two cases.

\begin{defn}
\label{defn:D}
The \mdfn{$D$-locus} is the space of all elements $\{a,b\}$ of $A_{n+1}$
with $a$ and $b$ in $\C_n^\perp$
such that 
\begin{enumerate}
\item
$a$ and $b$ are $\C$-orthogonal,
\item
$a$ and $\Ann(b)$ are orthogonal, and
\item
$b$ and $\Ann(a)$ are orthogonal.
\end{enumerate}
\end{defn}

\begin{remark}
\label{rem:D}
Since $\Ann(b)$ is a $\C$-subspace of $A_n$,
$a$ is orthogonal to $\Ann(b)$ if and only if $a$ is
$\C$-orthogonal to $\Ann(b)$.  
Similarly, $b$ is orthogonal to $\Ann(a)$ if and only if $b$
is $\C$-orthogonal to $\Ann(a)$.  Thus, conditions (2) and (3)
of Definition \ref{defn:D}
can be rewritten in terms of $\C$-orthogonality.

Also, $\Ann(b)^\perp$ is
equal to the image of left multiplication by $b$ 
(see Lemma \ref{lem:cancellation} below),
so condition (2) is also equivalent
to requiring that $a=bx$ for some $x$.  
Similarly, condition (3) is also equivalent
to requiring that $b=ay$ for some $y$.  
\end{remark}

The point of the following lemma is to determine precisely when
condition (4) of Proposition~\ref{prop:ann-intersect-H-perp} vanishes.

\begin{lemma}
\label{lem:D-locus-vanish}
Suppose that $a$ and $b$ belong to $\C_n^\perp$.
Then $\{a,b\}$ belongs to the $D$-locus if and only if
\[
(\beta^* - \alpha)\pi_{\C}(ab) + \pi_{\C}(ay-xb) = 0
\]
for all $\alpha$ and $\beta$ in $\C_n$, $x$ in $\Ann(a)$, and $y$ 
in $\Ann(b)$.
\end{lemma}

\begin{proof}
Since $\alpha$, $\beta$, $x$, and $y$ are independent,
the displayed expression vanishes 
if and only if
$\pi_{\C} (ab) = 0$, $\pi_{\C} (xb) = 0$ for
all $x$ in $\Ann a$, and $ay = 0$ 
for all $y$ in $\Ann b$.
The first equation just means that $a$ and $b$ are $\C$-orthogonal,
the second equation means that $b$ is $\C$-orthogonal to $\Ann(a)$,
and the third equation means that $a$ is $\C$-orthogonal to $\Ann(b)$.
\end{proof}

\begin{lemma}
\label{lem:D-locus-independent}
If $\{a,b\}$ is non-zero and does not belong to the $D$-locus, then
the dimension of $\Ann\{a,b\} \cap \HH_{n+1}^\perp$ is equal to
$\dim \Ann(a) + \dim \Ann(b)$.
\end{lemma}

\begin{proof}
Let $V$ be the subspace of $A_{n+1}$ consisting of all elements
of the form 
$\{ \alpha a + x, \beta b + y \}$, where $\alpha$ and $\beta$
belong to $\C_n$, $x$ belongs to $\Ann(a)$, and $y$ belongs to $\Ann(b)$.
The dimension of $V$ is
$\dim \Ann(a) + \dim \Ann(b) + 4$.
As in the proof of Theorem \ref{thm:ann-bracket-bound},
$V$ is a $\C_n$-vector space.

According to Proposition \ref{prop:ann-intersect-H-perp},
$\Ann\{a,b\} \cap \HH_{n+1}^\perp$ is contained in $V$.  In fact,
it is the subspace of $V$ defined by the two conjugate-linear equations
\begin{myequation}
\label{eq:Dli1}
\norm{a}^2 \alpha + \norm{b}^2 \beta^* = 0 \\
  \addtocounter{subsection}{1}
\label{eq:Dli2}
(\beta^* - \alpha)\pi_{\C}(ab) + \pi_{\C}(ay-xb) = 0.
\end{myequation}
Thus, we only need to show that Equations (\ref{eq:Dli1}) and (\ref{eq:Dli2})
are non-degenerate and independent.
Equation (\ref{eq:Dli1}) is non-degenerate 
because $\norm{a}$ or $\norm{b}$
is non-zero.  
Equation (\ref{eq:Dli2}) is non-degenerate 
by Lemma \ref{lem:D-locus-vanish}.

It remains to show that Equations (\ref{eq:Dli1}) and (\ref{eq:Dli2})
are independent.  There are three cases to consider, depending
on which part of Definition \ref{defn:D} fails to hold for $a$ and $b$.

If $a$ and $b$ are not $\C$-orthogonal, then $\pi_{\C}(ab)$ is non-zero.
Substitute the values $\alpha = -\norm{b}^2$, $\beta = \norm{a}^2$,
$x = 0$, and $y = 0$ into the two equations; note that
Equation (\ref{eq:Dli1}) is satisfied, while 
Equation (\ref{eq:Dli2}) is not satisfied because the left-hand side
equals $(\norm{a}^2 + \norm{b}^2) \pi_{\C}(ab)$.  This shows
that the two equations are independent because they have different solution
sets.

Next, suppose that $a$ is not orthogonal to $\Ann(b)$.
There exists an element $y_0$ of $\Ann(b)$ such that $a$ and $y_0$ are 
not $\C$-orthogonal.  This means that $\pi_{\C}(ay_0)$ is non-zero.
Substitute the values $\alpha = 0$, $\beta = 0$,
$x = 0$, and $y = y_0$ into the two equations; note that
Equation (\ref{eq:Dli1}) is satisfied, while 
Equation (\ref{eq:Dli2}) is not satisfied because the left-hand side
equals $\pi_{\C}(ay_0)$.  This shows
that the two equations are independent because they have different solution
sets.

Finally, suppose that $b$ is not orthogonal to $\Ann(a)$.
Similarly to the previous case, choose $x_0$ in $\Ann(a)$ such that
$\pi_{\C}(ax_0)$ is non-zero.
Substitute the values $\alpha = 0$, $\beta = 0$,
$x = x_0$, and $y = 0$ into the two equations; note that
Equation (\ref{eq:Dli1}) is satisfied, while 
Equation (\ref{eq:Dli2}) is not satisfied.
\end{proof}

\begin{thm}
\label{thm:ann-off-D-locus}
Let $a$ and $b$ be non-zero elements of $\C_n^\perp$.
If $\{a,b\}$ does not belong to the $D$-locus, 
then $\Ann\{a,b\}$ is contained in $\HH_{n+1}^\perp$.  
Moreover, the dimension of $\Ann\{a,b\}$ is
$\dim \Ann a +\dim \Ann b$.
\end{thm}

\begin{proof}
Recall from Lemma \ref{lem:zd-C-perp} 
that $\Ann\{a,b\}$ is a subspace of $\C_{n+1}^\perp$.
Also, $\HH_{n+1}^\perp$ is a codimension 2 subspace of $\C_{n+1}^\perp$.
Therefore, 
the codimension of 
$\Ann\{a,b\} \cap \HH_{n+1}^\perp$ in $\Ann\{a,b\}$ is at most 2.
Together with Lemma \ref{lem:D-locus-independent},
this implies that
the dimension of $\Ann\{a,b\}$ is at least
$\dim \Ann a + \dim \Ann b$ and at most 
$\dim \Ann a + \dim \Ann b + 2$.  However, the dimension
of $\Ann \{a,b\}$ is a multiple of 4 by Theorem \ref{thm:4dim},
so it must equal
$\dim \Ann a + \dim \Ann b$.  
This shows that $\Ann\{a,b\}$ equals $\Ann\{a,b\} \cap \HH_{n+1}^\perp$
because their dimensions are equal;
in other words, $\Ann\{a,b\}$ is contained in $\HH_{n+1}^\perp$.
\end{proof}

Theorem \ref{thm:ann-off-D-locus}
computes the dimension of $\Ann\{a,b\}$
for any $\{a,b\}$ that does not belong to the $D$-locus.  However,
it leaves something to be desired because it does not explicitly
describe $\Ann\{a,b\}$ as a subspace of $A_{n+1}$.  The difficulty arises
from our use of the fact that the dimension of $\Ann\{a,b\}$ is a multiple
of 4.

\begin{ques}
Describe $\Ann\{a,b\}$ explicitly when $\{a,b\}$ does not belong to
the $D$-locus.
\end{ques}

The rest of this section considers annihilators of elements that belong to
the $D$-locus.

\begin{lemma}
\label{lem:cancellation}
Suppose that $a$ and $b$ belong to $A_n$, and suppose that
$b$ is orthogonal to $\Ann(a)$.  There exists a unique element $x$
such that $ax = b$ and $x$ is orthogonal to $\Ann(b)$.
\end{lemma}

\begin{proof}
This is a restatement of Lemma \ref{lem:ann-im}.
\end{proof}

\begin{defn}
\label{defn:cancellation}
Let $a$ and $b$ belong to $A_n$, and suppose that $b$ is orthogonal
to $\Ann a$.  Then \mdfn{$\frac{b}{a}$} is the unique element
such that $a \frac{b}{a} = b$ and 
such that $\frac{b}{a}$ is orthogonal to $\Ann a$.
\end{defn}

Beware that the definition of $\frac{b}{a}$ is not symmetric.
In other words, it is not always true
that $\frac{b}{a} a = b$.

\begin{lemma}
\label{lem:frac-perp}
Let $a$ and $b$ be $\C$-orthogonal elements of $\C_n^\perp$, and suppose that 
$b$ is orthogonal to $\Ann(a)$.  
Then $\frac{b}{a}$ belongs to $\C_n^\perp$ and is 
$\C$-orthogonal to both $a$ and $b$.
\end{lemma}

\begin{proof}
If $a = 0$, then $\Ann a$ is all of $A_n$ so $b = 0$ and $\frac{b}{a}$
also equals $0$. In this case, the claim is
trivially satisfied.  Now assume that $a$ is non-zero.

For the first claim, note that 
$\langle a, a \frac{b}{a} \rangle_{\C} = \langle a, b \rangle_{\C} = 0$.
By Lemma \ref{lem:C-conj-linear}, this equals
$-\langle a^2, \frac{b}{a} \rangle^*_{\C}$. 
But $a^2$ is a non-zero real number, so $\frac{b}{a}$ 
is $\C$-orthogonal to $1$ as desired.

Next, note that $a \frac{b}{a} = b$ is 
orthogonal to $\C_n$, 
so $\left\langle a, \frac{b}{a} \right\rangle_{\C} = \pi_{\C}(b)$
is zero.
Also, compute that
\[
\textstyle\left\langle \frac{b}{a}, b \right\rangle_{\C} =
\left\langle \frac{b}{a}, a \frac{b}{a} \right\rangle_{\C} =
-\left\langle \left(\frac{b}{a}\right)^2, a \right\rangle^*_{\C}
\]
using Lemma \ref{lem:C-conj-linear}.
But $\left(\frac{b}{a}\right)^2$ is a real scalar,
which is $\C$-orthogonal to $a$ because we assumed that $a$ 
belongs to $\C_n^\perp$.
\end{proof}

\begin{thm}
\label{thm:ann-D-locus}
Let $a$ and $b$ be non-zero elements of $\C_n^\perp$, and 
suppose that $\{a,b\}$ belongs to the $D$-locus.  Then
$\Ann\{a,b\}$ is the $\C$-orthogonal direct sum of:
\begin{enumerate}
\item 
the space of all elements $\{x,y\}$ such that $x$ belongs to $\Ann(a)$
and $y$ belongs to $\Ann(b)$;
\item
the $\C$-linear span of the element
$\left\{ \norm{b}^2 a, -\norm{a}^2 b \right\}$;
\item
the $\C$-linear span of 
$\left\{ \frac{b}{a}, -\frac{a}{b} \right\} + \sqrt{2}(0,i_n)$,
where $\frac{b}{a}$ and $\frac{a}{b}$ are described
in Definition~\ref{defn:cancellation}.
\end{enumerate}
In particular, the dimension of $\Ann\{a,b\}$ is equal to
$\dim(\Ann a)+\dim(\Ann b)+4$.
\end{thm}

\begin{proof}
It follows from Proposition \ref{prop:ann-intersect-H-perp}
that $\Ann\{a,b\}$ contains the space described in part (1).
Recall that
Lemma \ref{lem:D-locus-vanish} implies that condition (4)
of Proposition \ref{prop:ann-intersect-H-perp} vanishes.

Next, note that $\left\{ \norm{b}^2 a, -\norm{a}^2 b \right\}$ satisfies the
conditions of Proposition \ref{prop:ann-intersect-H-perp}.
It corresponds to $\alpha = \norm{b}^2$, $\beta = -\norm{a}^2$,
$x = 0$, and $y = 0$.

Finally, we want to show that
$\{a,b\} \left\{\frac{b}{a},-\frac{a}{b} \right\}+ \sqrt{2}\{a,b\}(0,i_n)$
is zero.
Lemma \ref{lem:frac-perp} says that 
Proposition \ref{prop:bracket-multiply} applies to the first term,
which therefore equals 
$\sqrt{2} \left\{a \frac{b}{a}, -b\frac{a}{b} \right\}$.  This simplifies
to $\sqrt{2} \{ b, -a \}$.
Lemma \ref{lem:last-multiply} lets us compute that the second term is
$\sqrt{2} \{-b, a\}$, as desired.

We have now exhibited a subspace of $\Ann\{a,b\}$ whose dimension
is $\dim \Ann a + \dim \Ann b + 4$.  Theorem \ref{thm:ann-bracket-bound}
implies that we have described the entire annihilator.

Recall that Lemma \ref{lem:bracket-inner} describes how to compute
Hermitian inner products.  Using this lemma,
parts (1) and (2) are $\C$-orthogonal because
$a$ and $b$ are $\C$-orthogonal to $\Ann(a)$ and $\Ann(b)$ respectively.
Parts (1) and (3) are $\C$-orthogonal by 
Definition \ref{defn:cancellation}.
Parts (2) and (3) are $\C$-orthogonal by Lemma
\ref{lem:frac-perp}.
\end{proof}

%%%%%%%%%%%%%%%%%%%%%%%%%%%%%%%%%%%%%%%%%%%%%%%%%%%%%%%%%%%%%%%%%%

\section{The $D$-locus in $A_5$}
\label{sctn:D5}

The goal of this section is to explicitly understand the $D$-locus in $A_5$\
(see Definition \ref{defn:D}).  Unlike most of the rest of this paper,
this section uses computational techniques that apply in $A_4$ but have
not yet been made to work in general.

Let us consider whether elements of the form $\{a,0\}$ belong
to the $D$-locus.  If $a$ is non-zero, then part (2) of Definition \ref{defn:D}
fails.  Therefore, $\{a,0\}$ belongs to the $D$-locus only if $a = 0$.
Similarly, $\{0,b\}$ belongs to the $D$-locus only if $b = 0$.

From now on, we may suppose that $a$ and $b$ are non-zero.
If $b$ is not a zero-divisor, then it
is easy to determine whether $\{a,b\}$ belongs to the $D$-locus.
Namely, $b$ must be $\C$-orthogonal to $a$ and to $\Ann(a)$
because condition (2) of Definition \ref{defn:D} is vacuous.
By symmetry, a similar description applies when $a$ is not a zero-divisor.
Since annihilators in $A_4$ are well-understood 
\cite[Sections 11 and 12]{DDD}),
it is relatively straightforward to completely describe the elements
$\{a,b\}$ belonging to the $D$-locus in $A_5$ 
such that $a$ or $b$ is not a zero-divisor.

There is only one remaining case to consider.  
It consists of elements of the form $\{a,b\}$, where
$a$ and $b$ are both zero-divisors in $A_4$.  
We will focus on such elements in the rest of this section.
First we need some preliminary calculations in $A_4$.

\begin{lemma}
\label{lem:D5-1}
Let $a$ belong to $\C_3^\perp$.  If $b$ is $\C$-orthogonal
to $1$ and $a$, and $\alpha$ belongs to $\C_3$, then
the element $\{a , 0\}$ of $A_4$ is orthogonal to the annihilator
of $\{b, \alpha a\}$.
\end{lemma}

\begin{proof}
Let $c$ be the element of $A_3$ such that $bc = a$; in other words,
$c = -\frac{1}{\norm{b}^2} ba$.  Note that $c$ is $\C$-orthogonal to
both $a$ and $b$ by Lemma \ref{lem:ortho1}.

Using Proposition \ref{prop:bracket-multiply}, compute that
$\{b, \alpha a\} \{\frac{1}{\sqrt{2}} c, 0 \} = \{a, 0\}$.
Finally, use Lemma \ref{lem:ann-im} to conclude
that $\{a,0\}$ is orthogonal to $\Ann \{b, \alpha a\}$.
\end{proof}

\begin{lemma}
\label{lem:D5-2}
Let $a$ be a non-zero element of $\C_3^\perp$.
If a zero-divisor in $A_4$ is $\C$-orthogonal to $\{a,0\}$ 
and is orthogonal to $\Ann\{a,0\}$, then it is of the form
$\{b, \alpha a\}$, where $b$ is $\C$-orthogonal to $a$
and $\alpha$ belongs to $\C$.
\end{lemma}

\begin{proof}
Suppose that $x$ is a zero-divisor in $A_4$ that is $\C$-orthogonal
to $\{a,0\}$ and is orthogonal to $\Ann\{a,0\}$.
Write $x$ in the form 
$\{b, c\} + (\beta, \gamma)$,
where $b$ and $c$ belong to $\C_3^\perp$ 
while $\beta$ and $\gamma$ belong to $\C_3$.

Recall from \cite[Theorem 10.2]{DDD} that $\Ann\{a,0\}$ consists of 
elements of the form $\{0,y\}$, where $y$ is any element of $A_3$ that
is $\C$-orthogonal to $1$ and to $a$.
Since $x$ is orthogonal to $\Ann\{a,0\}$, Lemma \ref{lem:bracket-inner}
implies that $c$ is $\C$-orthogonal to all such $y$.
In other words,
$c$ belongs to the $\C$-linear span of $a$; i.e., $c = \alpha a$
for some $\alpha$ in $\C_3$.

Since $\{a,0\}$ and $x$ are $\C$-orthogonal, 
Lemma \ref{lem:bracket-inner} says that $b$ is $\C$-orthogonal to $a$.
Note, in particular, that $b$ and $c$ are $\C$-orthogonal.

Let $x_1 = b + c + \beta$
and $x_2 = -i_3 b + i_3 c + \gamma$ so that $x = (x_1, x_2)$.
Recall from \cite[Proposition~12.1]{DDD} that since $x$ is a zero-divisor,
$x_1$ and $x_2$ are 
imaginary orthogonal elements of $A_3$
with the same norm.  

Multiplication by $i_3$ preserves norms in $A_3$.  
Since $x_1$ and $x_2$ have the same norm, it follows that $\beta$
and $\gamma$ have the same norm.  This uses that 
$b$ and $c$ are orthogonal, as we have already shown.

Since $x_1$ and $x_2$ are orthogonal, it follows that
$\beta$ and $\gamma$ are orthogonal.  This uses that
$b$ and $c$ are each orthogonal to both $i_3 b$ and $i_3 c$
since $b$ and $c$ are $\C$-orthogonal.

Next, since $x_1$ and $x_2$ are imaginary, it follows that
$\beta$ and $\gamma$ are $\R$-scalar multiples of $i_3$.
We have shown that $\beta$ and $\gamma$ are both orthogonal and parallel
and also have the same norm.
It follows that $\beta$ and $\gamma$ are both zero.
\end{proof}

\begin{prop}
\label{prop:D5}
Let $a$, $b$, and $c$ belong to $\C_3^\perp$, and suppose that $a$ is
non-zero. Suppose also that $\{b,c\}$ is a zero-divisor in $A_4$.
The element
$\{ \{a,0\}, \{b,c\} \}$ belongs to the $D$-locus in $A_5$
if and only if $b$ is $\C$-orthogonal to $a$ and $c$ belongs to
the  $\C$-linear span of $a$.
\end{prop}

\begin{proof}
First suppose that $b$ is $\C$-orthogonal to $a$ and $c$ belongs to 
the $\C$-linear span of $a$.  
Lemma \ref{lem:bracket-inner} implies that $\{a,0\}$ and $\{b,c\}$ are
$\C$-orthogonal.  

By \cite[Theorem 10.2]{DDD}, 
the annihilator of $\{a,0\}$ consists of elements
of the form $\{0,y\}$, where $y$ is $\C$-orthogonal to $1$ and $a$.
Therefore, Lemma \ref{lem:bracket-inner} implies that $\{b,c\}$ is
orthogonal to $\Ann \{a,0\}$.

Lemma \ref{lem:D5-1} implies that $\{a,0\}$ is orthogonal to 
$\Ann \{b,c\}$.  This finishes one implication.

For the other implication, suppose that 
$\{ \{a,0\}, \{b,c\} \}$ belongs to the $D$-locus in $A_5$.
Lemma \ref{lem:D5-2} implies that $b$ is $\C$-orthogonal to $a$
and that $c$ belongs to the $\C$-linear span of $a$.
\end{proof}

Suppose that $a = (a_1, a_2)$ is a zero-divisor in $A_4$.
We recall from \cite[Sections 11 and 12]{DDD} some algebraic properties of $a$.
First of all, $a_1$ and $a_2$ are imaginary orthogonal elements of $A_3$
with the same norm.  
The $\R$-linear span of $1$, $a_1$, $a_2$, and $a_1 a_2$
is a 4-dimensional subalgebra $\llangle a_1, a_2 \rrangle$ 
of $A_3$ that is isomorphic to the quaternions.
The notation indicates that the subalgebra is generated by $a_1$ and $a_2$.

The annihilator $\Ann(a)$ is a four-dimensional subspace of $A_4$
consisting of all elements of the form $(y, -cy)$, where 
$c$ is the fixed unit vector with the same direction as $a_1 a_2$
and $x$ ranges over 
the orthogonal complement of $\llangle a_1, a_2 \rrangle$.
The subspace $\llangle a_1, a_2 \rrangle \times \llangle a_1, a_2 \rrangle$
is orthogonal to $\Ann(a)$.
Let \mdfn{$\Eig_2(a)$} be the orthogonal complement of $\Ann(a)$ and
$\llangle a_1, a_2 \rrangle \times \llangle a_1, a_2 \rrangle$.
This space
consists of all elements of the form $(y,cy)$,
where $c$ and $x$ are as above.  Direct calculation shows
that $\Eig_2(a)$ is equal to the space
of all elements $b$ of $A_4$ such that $a(ab) = -2b$.
From this perspective, it is the 2-eigenspace of the composition
of left multiplication by $a$ and left multiplication by $a^* = -a$.

\begin{cor}
\label{cor:D5}
Let $a = (a_1, a_2)$ and $b = (b_1, b_2)$ be zero-divisors in $A_4$.
Then $\{a,b\}$ belongs to the $D$-locus in $A_5$ if and only if
$b$ belongs to the 
$\R$-linear span of 
$(a_1, -a_2)$, $(a_2, a_1)$, and $\Eig_2(a)$.
\end{cor}

\begin{proof}
Since $a_1$ and $a_2$ are orthogonal and have the same norm,
there exists an imaginary element $c$ of unit length such that
$a_2 = ca_1$.  There exists an automorphism of $A_3$
that takes $c$ to $-i_3$.  Therefore, we may assume that $c = -i_3$.
In other words, we may assume that $a = \{a_1, 0\}$.

Then $\Eig_2(a)$ is equal to the space of all elements of the
form $\{y,0\}$, where $y$ is $\C$-orthogonal to $1$ and $a$.
Also, $\{0,a\}$ equals $(a_1, -a_2)$, so the $\C_4$-linear span
of $\{0, a\}$ is the same as the $\R$-linear span of 
$(a_1, -a_2)$ and $(a_2, a_1)$.

Finally, apply Proposition \ref{prop:D5}.
\end{proof}

Recall that $V_2(\R^7)$ is the space of orthonormal 2-frames
in $\R^7$.  In the following theorem, we identify this space
with the space of 
elements $(a_1, a_2)$ of $A_4$ such that $a_1$ and $a_2$ are
orthogonal imaginary unit vectors in $A_3$.

\begin{thm}
\label{thm:D5}
Consider the space $X$ consisting of all elements $\{a,b\}$ belonging
to the $D$-locus in $A_5$
such that $a$ and $b$ are zero-divisors with unit length.
Let $\xi$ be the 4-plane bundle over $V_2(\R^7)$ whose
unit sphere bundle has total space diffeomorphic to 
the 14-dimensional compact simply connected Lie group $G_2$
(see \cite[Section 7]{DDD}).
Then $X$ is diffeomorphic to the unit sphere bundle of $\xi \oplus 2$,
where $\xi \oplus 2$ is the fiberwise sum of the vector bundle $\xi$ with
the trivial 2-dimensional bundle.
\end{thm}

\begin{proof}
First, identify $V_2(\R^7)$ with the space of all zero-divisors
in $A_4$ with unit length.  Let $\eta$ be the bundle over
$V_2(\R^7)$ whose fiber over $a$ is 
the space of all ordered pairs $(a,b)$ such that
$b$ is a unit length element of $\Eig_2(a)$.
The bundle $\xi$ is also a bundle over $V_2(\R^7)$, but the 
fiber over $a$ is
the space of all ordered pairs $(a,b)$ such that
$b$ is a unit length element of $\Ann(a)$.

Using the notation in the paragraphs preceding Corollary \ref{cor:D5},
the isomorphism
$\Eig_2(a) \map \Ann(a) \colon (y,cy) \mapsto (y,-cy)$ 
induces an isomorphism from $\eta$
to $\xi$.

Next consider the space of all ordered pairs $(a,b)$ such that
$a$ is a unit length zero-divisor and $b$ belongs to the 
$\R$-span of $(a_1, -a_2)$ and $(a_2, a_1)$, where $a = (a_1, a_2)$.
The map that takes $(a,b)$ to $a$ is a trivial
2-plane bundle.

Corollary \ref{cor:D5} shows that $X$ is the unit sphere bundle
of $\eta \oplus 2$.  
\end{proof}

\begin{remark}
An obvious consequence of Theorem \ref{thm:D5} is that $X$ is 
diffeomorphic to the total space of an $S^5$-bundle over
$V_2(\R^7)$.  This bundle is the fiberwise double suspension of
the usual $S^3$-bundle over $V_2(\R^7)$ that is used to construct 
$G_2$.
\end{remark}

%%%%%%%%%%%%%%%%%%%%%%%%%%%%%%%%%%%%%%%%%%%%%%%%%%%%%

\section{Stability}
\label{sctn:stable}

Sections \ref{sctn:ann-H-perp} and \ref{sctn:D-locus} described many
properties of annihilators of elements of the form $\{a,b\}$.  
This section exploits these properties to study large annihilators,
i.e., annihilators in $A_n$ whose dimension is at least $2^{n-1}$.

We begin with a result that could have been included in \cite{DDD},
but its significance was not apparent at the time.

\begin{thm}
\label{thm:top-half}
Let $n \geq 3$, and let $a$ belong to $A_n$.  
If the dimension of $\Ann(a)$ is at 
least $2^{n-1}$, then 
$a$ belongs to $\HH_n^\perp$.
\end{thm}

\begin{proof}
Let $a = (b,c)$.
We claim that $b$ and $c$ are both zero-divisors; otherwise,
\cite[Lemma 9.9]{DDD} would imply that $\Ann(a)$ has dimension at most 
$2^{n-1} - 1$.
Lemma \ref{lem:zd-C-perp} implies that $b$ and $c$ belong to $\C_{n-1}^\perp$.
\end{proof}

Theorem \ref{thm:top-half} is important in the following way.
When searching for zero-divisors with large annihilators, i.e.,
with annihilators whose dimension is at least half the dimension of $A_n$,
one need only look in $\HH_n^\perp$.  Fortunately, 
Sections \ref{sctn:ann-H-perp} and \ref{sctn:D-locus} study
zero-divisors in $\HH_n^\perp$ in great detail.

Next we show by construction that the bound of Theorem \ref{thm:top-half}
is sharp in the sense that there exist elements of $A_n$ that do
not belong to $\HH_n^\perp$ but whose annihilators have
dimension $2^{n-1} - 4$.
Recall that an element $a$ of $A_n$ is alternative if
$a \cdot ax = a^2 x$ for all $x$.  For every $n$, there exist
elements of $A_n$ that are alternative.  
For example, a straightforward computation
shows that if $a$ is an alternative element of $A_{n-1}$,
then $(a,0)$ is an alternative element of $A_n$.

\begin{prop}
\label{prop:Dugger-ex}
Let $a$ be any non-zero alternative element of $\C_{n-1}^\perp$
such that $\norm{a} = 1$.
Then $\Ann(i_{n-1}, a)$ is equal to the set of all
elements of the form $(x, ai_{n-1} \cdot x)$
such that $x$ is $\C$-orthogonal to $1$ and to $a$.
In particular, the dimension of $\Ann(i_{n-1}, a)$ is equal
to $2^{n-1} - 4$.
\end{prop}

\begin{proof}
Let $x$ be $\C$-orthogonal to both $1$ and $a$.
Using Lemma \ref{lem:C-conj-linear},
compute that the product
$(i_{n-1}, a) (x, ai_{n-1} \cdot x )$ is always zero.
We have exhibited a subspace of $\Ann(i_{n-1},a)$ that has
dimension $2^{n-1} - 4$.
By Theorem \ref{thm:top-half}, this subspace must be equal to
$\Ann(i_{n-1},a)$.
\end{proof}

A proof of Proposition \ref{prop:Dugger-ex} also appears in
\cite[Theorem 4.4]{M3}.

\begin{ques}
Find all of the elements of $A_n$ that have annihilators of
dimension $2^{n-1} - 4$.
\end{ques}

The paper \cite{DDD} began an exploration of the largest annihilators
in $A_n$.  
Recall from Theorem \ref{thm:4dim} that
the annihilators in $A_n$ have dimension
at most $2^n - 4n +4$.  Moreover, Theorem 15.7 of \cite{DDD} gives a
complete description of the elements 
whose annihilators have dimension equal to this upper bound.
The rest of this section provides more results in a similar vein.

\begin{defn}
\label{defn:T}
Let $n \geq 4$, and let $c$ be a multiple of $4$
such that $0 \leq c \leq 2^n - 4n$.
The space \mdfn{$T^c_n$} is the space of elements of length one in $A_n$
whose annihilators have dimension at least
$(2^n - 4n + 4) - c$.  
\end{defn}

This is a change in the definition of $T^c_n$ from that
used in \cite{DDD}.  
The elements of $T^c_n$ are unit length zero-divisors 
whose annihilators
are within $c$ dimensions of the largest possible value.
The space $T^0_4$ is diffeomorphic to the Stiefel manifold
$V_2(\R^7)$ of orthonormal 2-frames in $\R^7$
\cite[Section 12]{DDD}.

We have imposed the condition $n \geq 4$ in order to avoid trivial exceptions
to our results involving well-known properties of $A_n$ for $n \leq 3$.
Also, we have imposed the condition $c \leq 2^n - 4n$ to ensure that
every element of $T_n^c$ is always a zero-divisor.

It follows from Lemma \ref{lem:zd-C-perp} that
$T_n^c$ is contained in $\C_n^\perp$.  Thus, if $a$ and $b$
lie in $T_n^c$, then it makes sense to talk about $\{a,b\}.$
Note that if $a$ is in $T_n^c$ then $\{a,0\}$ and $\{0,a\}$ lie in
$T_{n+1}^c.$  This is because, according to Proposition
\ref{prop:ann-bracket-special},
both $\Ann\{a,0\}$ and $\Ann\{0,a\}$ have dimension
equal to $\dim \Ann(a) + 2^n - 4$.  Consequently, $T_{n+1}^c$ contains
a disjoint union of two copies of $T_n^c.$

\begin{defn}
The space $T_n^c$ is \mdfn{stable} if $T_{n+1}^c$ is diffeomorphic to the
space of elements of the form $\{a,0\}$ or $\{0,a\}$ such that
$a$ belongs to $T_n^c$.
\end{defn}

For $n \geq 4$, the space $T_n^0$ is stable \cite[Proposition 15.6]{DDD};
a vastly simpler proof appears below.
In fact, our goal is to completely determine which spaces $T_n^c$ are stable.

\begin{prop}
\label{prop:stability}
Let $n \geq 4$, and let $a$ belong to $A_n$.  
If the dimension of $\Ann(a)$ is at least $2^n - 8n + 24$,
then $a$ is of the form $\{b,0\}$ or $\{0,b\}$ with
$b$ in $\C_{n-1}^\perp$.
\end{prop}

\begin{proof}
Suppose that the dimension of $\Ann(a)$ is at least $2^n - 8n + 24$.
Note that $2^{n-1} \leq 2^n - 8n + 24$,
so the dimension of $\Ann(a)$ is at least $2^{n-1}$.
By Theorem \ref{thm:top-half},
$a$ belongs to $\HH_n^\perp$.

Write $a =
\{x,y\}$ for some $x$ and $y$ in $\C_{n-1}^\perp$.  Assume for
contradiction that both $x$ and $y$ are non-zero.  By
Theorem~\ref{thm:ann-bracket-bound}, the dimension of $\Ann(a)$ is at
most $\dim \Ann(x) + \dim \Ann(y) + 4$.  But the dimensions of
$\Ann(x)$ and $\Ann(y)$ are at most $2^{n-1} - 4n + 8$ by 
Theorem \ref{thm:4dim}, so the dimension of $\Ann(a)$
is at most $2^n - 8n + 20$.  This is a contradiction, so either
$x$ or $y$ is zero.
\end{proof}

\begin{prop}
\label{prop:stability2}
If $n \geq 4$, $c \geq 0$, and $n \geq \frac{c}{4}+4$, then $T^c_n$ is stable.
\end{prop}

\begin{proof}
It follows from the inequalities that $c \leq 2^n - 4n$.

Let $a$ belong to $T^c_{n+1}$.
Note that $2^{n+1} - 4n - c \geq 2^{n+1} - 8n + 16$, 
so $\Ann(a)$ has dimension at least $2^{n+1} - 8(n+1) + 24$.
Proposition \ref{prop:stability} implies that $a$ is of the form
$\{b,0\}$ or $\{0,b\}$.
The result then follows directly from Proposition
\ref{prop:ann-bracket-special}.
\end{proof}

\begin{lemma}
\label{lem:top-dim-D-locus}
For $n \geq 3$, there exists an element $\{a,b\}$ belonging to
the $D$-locus in $A_{n+1}$
such that $a$ and $b$ are elements of $\C_{n}^\perp$ whose annihilators
have dimension $2^{n} - 4n + 4$.
\end{lemma}

\begin{proof}
The proof is by induction on $n$.  The base case is $n = 3$.
Since every element of $A_3$ has a trivial annihilator,
this case just requires us to choose two $\C$-orthogonal elements
from the 6-dimensional space $\C_3^\perp$.

Now suppose that $a'$ and $b'$ are elements of $\C_n^\perp$
whose annihilators
have dimension $2^n - 4n + 4$.  Suppose also that
$\{a', b'\}$ belongs to the $D$-locus in $A_{n+1}$.

Consider the elements $a = \{a',0\}$ and $b =\{b',0\}$ of $A_{n+1}$.
By Proposition \ref{prop:ann-bracket-special} and the induction
assumption, $a$ and $b$ have annihilators of dimension
$2^{n+1} - 4(n+1) + 4$, as desired.  

It remains to show that $\{a,b\}$ belongs to the $D$-locus in $A_{n+2}$.
By Lemma \ref{lem:bracket-inner} and the induction assumption,
$a$ and $b$ are $\C$-orthogonal.  
Proposition \ref{prop:ann-bracket-special} describes
$\Ann(a)$ and $\Ann(b)$.  By inspection of this description,
$b$ is $\C$-orthogonal to $\Ann(a)$ because $b'$ is $\C$-orthogonal
to $\Ann(a')$ by the induction assumption.
Similarly, $a$ is $\C$-orthogonal to $\Ann(b)$.
\end{proof}

\begin{lemma}
\label{lem:top-dim-D-locus2}
For $n \geq 4$, there exist non-zero elements $a$ and $b$ in $\C_n^\perp$
such that $\Ann\{a,b\}$ has dimension $2^{n+1} - 8n + 12$.
\end{lemma}

\begin{proof}
By Lemma \ref{lem:top-dim-D-locus}, there exist non-zero elements of
$\C_n^\perp$ such that $\{a,b\}$ belongs to the $D$-locus in $A_{n+1}$
and such that $\Ann(a)$ and $\Ann(b)$ both have dimension
$2^n - 4n + 4$.
Now apply Theorem \ref{thm:ann-D-locus} to conclude that
$\Ann\{a,b\}$ has dimension $2^{n+1} - 8n + 12$.
\end{proof}

\begin{remark}
Lemma \ref{lem:top-dim-D-locus2} shows that the bound of 
Proposition \ref{prop:stability} is sharp.  Substitute $n-1$ for $n$
in the lemma to construct an element of $A_n$ whose annihilator
has dimension $2^n - 8n + 20$.
\end{remark}

\begin{prop}
\label{prop:not-stable}
Let $n \geq 4$, let $c \leq 2^n - 4n$,
and let $n \leq \frac{c}{4} + 3$.  Then $T_n^c$ is not stable.
\end{prop}

\begin{proof}
Note that $2^{n+1} - 4(n+1) + 4 - c \leq 2^{n+1} - 8n + 12$.
Now apply Lemma \ref{lem:top-dim-D-locus2} to construct an element
$\{a,b\}$ belonging to $T_{n+1}^c$ such that both $a$ and $b$ are non-zero.
\end{proof}

\begin{thm}
\label{thm:stable-dim}
Let $n \geq 4$, and let $c$ be a multiple of $4$
such that $0 \leq c \leq 2^n - 4n$.  Then $T_n^c$ is stable if and only 
if $n \geq \frac{c}{4} + 4$.
\end{thm}

\begin{proof}
Combine Propositions \ref{prop:stability2} and \ref{prop:not-stable}.
\end{proof}

We give two illustrations of the theorem.

\begin{cor}
\label{cor:fuck-matrices}
The space of zero-divisors in $A_5$ whose annihilators are 16-dimensional
is diffeomorphic to two disjoint copies of $V_2(\R^7)$.
\end{cor}

\begin{proof}
Apply Theorem \ref{thm:stable-dim} with 
$n = 5$ and $c = 0$.  
\end{proof}

Corollary \ref{cor:fuck-matrices} is the same as
\cite[Corollary 14.7]{DDD}.
The proof is vastly more graceful than the one in \cite{DDD}.
This demonstrates the power of our computational perspective.

\begin{cor}
The space of zero-divisors in $A_6$ whose annihilators are 
at least 40-dimensional
is diffeomorphic to two disjoint copies of the space of
zero-divisors in $A_5$ whose annihilators are at least 12-dimensional.
\end{cor}

%%%%%%%%%%%%%%%%%%%%%%%%%%%%%%%%%%%%%%%%%%%%%%%%%%%%%%%%%%%%%%%%5

\bibliographystyle{amsalpha}

\end{document}